\input amstex
\documentstyle{amsppt}
\topmatter
\magnification=\magstep1
\pagewidth{5.2 in}
\pageheight{6.7 in}
\abovedisplayskip=10pt
\belowdisplayskip=10pt
\NoBlackBoxes
\title
Some identities for the Bernoulli, the Euler and the Genocchi numbers and polynomials
\endtitle
\author  Taekyun Kim  \endauthor
\affil\rm{{Division of General Education-Mathematics,}\\
{ Kwangwoon University, Seoul 139-701, S. Korea}\\
{e-mail: tkkim$\@$kw.ac.kr}}
\endaffil

\abstract{ The purpose of this paper is to give  some new identities for the Bernoulli,
 the Euler and the Genocchi numbers and polynomials.
}
\endabstract
\thanks 2000 Mathematics Subject Classification  11S80, 11B68 \endthanks
\thanks Key words and phrases: Euler number, $p$-adic invariant integrals, zeta function, $p$-adic fermionic integrals \endthanks
\rightheadtext{   } \leftheadtext{ Euler numbers and
polynomials }
\endtopmatter

\document

\head 1. Introduction \endhead
Let $p$ be a fixed odd prime and let  $\Bbb Z_p$, $\Bbb Q_p$, $\Bbb C$ and $\Bbb C_p$ denote the ring of $p$-adic rational integers,
the field of $p$-adic rational numbers, the complex number field and the completion of algebraic closure of $\Bbb Q_p $.
The $p$-adic absolute value in $\Bbb C_p$ is normalized so that
$|p|_p=\frac1p$.
For $ f \in UD( \Bbb Z_p) $, let us start with the expression
$$ \sum_{ 0 \leq j < p^N } (-1)^j f(j) = \sum_{0 \leq j < p^N} f(j)\mu ( j + p^N \Bbb Z_p )$$
representing analogue of Riemann's sums for $f$, cf.[1-10].

The fermionic $p$-adic invariant integral of $f$ on $\Bbb Z_p$ will be
defined as the limit $( N \rightarrow \infty )$ of these sums,
which it exists. The fermionic $p$-adic invariant integral of a function
$ f \in UD(\Bbb Z_p)$ is defined in [1, 4, 5, 10] as follows:
$$I(f)= \int_{\Bbb Z_p} f(x) d \mu(x)=\lim_{N\rightarrow \infty}\sum_{0 \leq j < p^N} f(j)\mu ( j + p^N \Bbb Z_p ) = \lim_{N \rightarrow \infty}
\sum_{0 \leq j < p^N } f(j) (-1)^j .\tag1$$
From (1), we note that
$$I(f_1) + I(f)=2f(0), \text{ where $f_1(x)=f(x+1)$}. $$
By using integral iterative method, we also easily see that
$$I(f_n)+(-1)^{n-1}I(f)=2\sum_{l=0}^{n-1}(-1)^{n-1-l}f(l), \text{ where $f_n(x)=f(x+n)$ for  $n \in \Bbb N$}. \tag2$$
For
$d$ a fixed positive even integer with $(p,d)=1$, let
$$X=X_d=\varprojlim_N \Bbb Z/dp^N\Bbb Z , \;\;X_1=\Bbb Z_p,$$
$$X^*=\bigcup\Sb 0<a<dp\\ (a,p)=1\endSb a+dp\Bbb Z_p,$$
$$a+dp^N\Bbb Z_p=\{x\in X\mid x\equiv a\pmod{dp^N}\},$$
where $a\in \Bbb Z$ lies in $0\leq a<dp^N$, (see [1-10]).
Recently, several authors have studied the properties of $p$-adic invariant integral on $\Bbb Z_p$ related to Euler numbers
and polynomials (see [1, 2, 4, 5, 7, 8, 9, 10]). The integral equation (2) for $n \equiv 1$ $(mod \ 2)$ are useful to study the congruence
of Euler numbers and polynomials (see [1, 2, 4, 5, 7, 8, 9, 10]).
In this paper, we consider the equations of the fermionic $p$-adic invariant integral on $\Bbb Z_p$ for $n\equiv 0$ $(mod \   2)$.
From those equations of the fermionic $p$-adic invariant integral on $\Bbb Z_p$ for $n\equiv 0$ $(mod  \  2)$, we derive some interesting and valuable identities for the Euler, the Genocchi and the Bernoulli numbers and polynomials.

\head 2. some identities of  the Bernoulli, the Euler and the Genocchi numbers and polynomials
\endhead

If $n\equiv 0$ $(mod  \  2)$ in (2), then we see that
$$I(f_n)-I(f)=2\sum_{l=0}^{n-1}(-1)^{l-1}f(l), \text{ $f_n(x)=f(x+n)$}. \tag3$$
Let us take $f(x)=e^{tx}$. Then we have
$$\int_{\Bbb Z_p}e^{xt}d\mu(x)
=\frac{2 \sum_{l=0}^{d-1}(-1)^{l-1}e^{lt}}{e^{dt}-1}, \text{ for $d\in\Bbb N$ with $d\equiv 0$ $(mod   \   2)$}. \tag4$$
It is easy to see that
$$\frac{2 \sum_{l=0}^{d-1}(-1)^{l-1}e^{lt}}{e^{dt}-1}=\frac{2}{e^t +1}=\sum_{l=0}^{\infty}E_n \frac{t^n}{n!}, \tag5 $$
where $E_n$ are the $n$-th Euler numbers.

From (4), we can also derive
$$\aligned
&\int_{\Bbb Z_p}e^{xt}d\mu(x)=\frac{2\sum_{l=0}^{d-1}(-1)^{l-1}e^{lt}}{e^{dt}-1}=2\sum_{l=0}^{d-1}\left(\frac{dt}{e^{dt}-1}\right)\frac{e^{lt}}{dt}\\
&= \frac{1}{dt}2\sum_{l=0}^{d-1}(-1)^{l-1}\left(\sum_{n=1}^{\infty}B_n(\frac{l}{d})d^n \frac{t^n}{n!}\right)
=\sum_{n=0}^{\infty}\left(2d^n\sum_{l=0}^{d-1}(-1)^{l-1}\frac{B_{n+1}(\frac{l}{d})}{n+1} \right)\frac{t^n}{n!},
\endaligned\tag6$$
where $B_n(x)$ are the $n$-th Bernoulli polynomials.

Thus, we have
$$\int_{\Bbb Z_p}x^n d\mu(x)=\frac{2d^n}{n+1}\sum_{l=0}^{d-1}(-1)^{l-1}B_{n+1}(\frac{l}{d}),
 \text{ where $d\in\Bbb N$ with $d\equiv 0$ $(mod  \   2)$}.
\tag7$$
By (4) and (5), we also see that
$$ \int_{\Bbb Z_p}e^{tx} d \mu(x)=\frac{2}{e^t +1}=\sum_{n=0}^{\infty} E_n \frac{t^n}{n!}, \tag8$$
where $E_n$ are the $n$-th Euler numbers.

From (7) and (8), we obtain the following theorem.
\proclaim{Theorem 1} For $n\in \Bbb Z_{+},$  $d \in \Bbb N$ with $ d\equiv 0$ $(mod   \   2)$,  we have
$$ \frac{E_n}{2}=\frac{d^n}{n+1}\sum_{l=0}^{d-1}(-1)^{l-1}B_{n+1}(\frac{l}{d}).$$
\endproclaim
It is easy to show that
$$\aligned
&\int_{\Bbb Z_p}e^{(x+y)t}d\mu(y)=\frac{2\sum_{l=0}^{d-1}(-1)^{l-1}e^{lt}}{e^{dt}-1}e^{xt}=\frac{2}{dt}\sum_{l=0}^{d-1}(-1)^{l-1}
\left(\frac{dt}{e^{dt}-1}\right)e^{(l+x)t}\\
&=2\sum_{l=0}^{d-1}(-1)^{l-1}\sum_{n=0}^{\infty}\frac{B_{n+1}(\frac{l+x}{d})}{n+1}\frac{d^nt^n}{n!}
=\sum_{n=0}^{\infty}\left(\frac{2d^n}{n+1}\sum_{l=0}^{d-1}(-1)^{l-1}B_{n+1}(\frac{l+x}{d})\right)\frac{t^n}{n!}.
\endaligned \tag9$$
From (3), we can easily derive the following equation (3).
$$\int_{\Bbb Z_p}e^{(x+y)t}d\mu(y)=\frac{2\sum_{l=0}^{d-1}(-1)^{l-1}e^{lt}}{e^{dt}-1}e^{xt}=\frac{2}{e^{t}+1}e^{xt}=\sum_{n=0}^{\infty}E_n(x)\frac{t^n}{n!}.
\tag10$$
By (3), (9) and (10), we have
$$\frac{E_n(x)}{2}=\frac{1}{2}\int_{\Bbb Z_p}(y+x)^nd\mu(y)=\frac{d^n}{n+1}\sum_{l=0}^{d-1}(-1)^{l-1}B_{n+1}(\frac{l+x}{d}), $$
 and
 $$ \frac{1}{2}\left(E_n(d)-E_n \right) = \frac{d^n}{n+1}
 \left(\sum_{l=0}^{d-1}(-1)^{l-1}\left(B_{n+1}(\frac{l}{d}+1)-B_{n+1}(\frac{l}{d})\right)\right)=\sum_{l=0}^{d-1}(-1)^{l-1}l^n.$$
 Therefore, we obtain the following theorem.
 \proclaim{ Theorem 2} For $d \in \Bbb N $ with $ d\equiv 0$ $(mod  \   2)$,  $ n\in \Bbb Z_{+}$, we have
 $$\sum_{l=0}^{d-1}(-1)^{l-1}\left(\frac{l}{d}\right)^n
 =\frac{1}{n+1}\left( \sum_{l=0}^{d-1}(-1)^{l-1}\left(B_{n+1}(\frac{l}{d}+1)-B_{n+1}(\frac{l}{d})\right)\right), $$
 and
 $$ \frac{1}{2}\left( E_n(d)-E_n \right)=\sum_{l=0}^{d-1}(-1)^{l-1}l^n.$$
    \endproclaim
Let us define  the $n$-th Genocchi polynomials as follows:
for $d\in \Bbb N$ with $d\equiv 0$ $(mod  \   2)$,
$$t\int_{\Bbb Z_p}e^{(x+y)t}d\mu(y)=\frac{2t\sum_{l=0}^{d-1}(-1)^{l-1}e^{lt}}{e^{dt}-1}e^{xt}=\sum_{n=0}^{\infty}G_n(x)\frac{t^n}{n!}.\tag11$$
Thus, we note that
$$\aligned
\frac{2t\sum_{l=0}^{d-1}(-1)^{l-1}e^{lt}}{e^{dt}-1}e^{xt}&=\frac{2}{d}\frac{dt\sum_{l=0}^{d-1}(-1)^{l-1}e^{lt}}{e^{dt}-1}e^{xt}
=\frac{2}{d}\sum_{l=0}^{d-1}(-1)^{l-1}\sum_{n=0}^{\infty}B_n(\frac{l+x}{d})\frac{t^nd^n}{n!}\\
&=\sum_{n=0}^{\infty}\left(2d^{n-1}\sum_{l=0}^{d-1}(-1)^{l-1}B_n(\frac{l+x}{d})\right)\frac{t^n}{n!}.
\endaligned\tag12$$
By comparing coefficients on the both sides of (11) and (12), we obtain the following theorem.

\proclaim{ Theorem 3} For $d\in\Bbb N$ with $d\equiv 0$ $(mod  \   2)$, we have
$$\frac{G_n(x)}{2}= d^{n-1}\sum_{l=0}^{d-1}(-1)^{l-1}B_n(\frac{l+x}{d}), \text{ for $n\in\Bbb Z_{+}$}.$$
\endproclaim
Recently, several authors have studied the generalized Euler numbers and polynomials attached to the Dirichlet's character with
odd conductor (see [1, 2, 4, 5, 7, 8, 9, 10, 11, 12]). Now, we consider the generalized Euler polynomials attached to the Dirichlet's
 character with even conductor. For $d \in \Bbb N$ with $d\equiv 0$ $(mod  \   2)$, let $\chi$ be the Dirichlet's character with
 conductor $d$.
From (2), we note that
$$\int_{X}\chi(y)e^{(x+y)t}d\mu(y)=\left(\frac{2\sum_{l=0}^{d-1}(-1)^{l-1}\chi(l)e^{lt}}{e^{dt}-1}\right)e^{xt}. \tag13$$
Let us define the generalized Euler polynomials for the Dirichlet's character with even conductor as follows: for $d\in\Bbb N$ with
$d\equiv 0$ $(mod  \   2)$,
$$\left(\frac{2\sum_{l=0}^{d-1}(-1)^{l-1}\chi(l)e^{lt}}{e^{dt}-1}\right)e^{xt}=\sum_{n=0}^{\infty}E_{n, \chi}(x)\frac{t^n}{n!}.\tag14$$
From (13) and (14), we can also derive the following equation:
$$\int_{X}\chi(y)(x+y)^n d\mu(y)=E_{n, \chi}(x), \text{ for $n\in \Bbb Z_{+}$}.$$
The $n$-th Genocchi polynomials are also defined by
$$\left(\frac{2t\sum_{l=0}^{d-1}(-1)^{l-1}\chi(l)e^{lt}}{e^{dt}-1}\right)e^{xt}=\sum_{n=0}^{\infty}G_{n, \chi}(x)\frac{t^n}{n!},
\text{ where $d\in\Bbb N$ with $d\equiv 0$ $(mod  \  2)$}.\tag15$$
By (15), we see that
$$\aligned
\sum_{n=0}^{\infty}G_{n,\chi}(x)\frac{t^n}{n!}&=\frac{2}{d}\sum_{l=0}^{d-1}(-1)^{l-1}\chi(l)\frac{dt}{e^{dt}-1}e^{(l+x)t}
=\frac{2}{d}\sum_{l=0}^{d-1}(-1)^{l-1}\chi(l)\sum_{n=0}^{\infty}d^n B_n(\frac{l+x}{d})\frac{t^n}{n!}\\
&=\sum_{n=0}^{\infty} \left\{ 2d^{n-1}\sum_{l=0}^{d-1}(-1)^{l-1}\chi(l)  B_n(\frac{l+x}{d}) \right\}\frac{t^n}{n!}.
\endaligned\tag16$$
By comparing coefficients on the both sides of (16), we obtain the following theorem.
\proclaim{ Theorem 4}
Let $d\equiv 0$ $(mod  \  2)$. Then we have
$$\frac{G_{n, \chi}(x)}{2}=d^{n-1}\sum_{l=0}^{d-1}(-1)^{l-1}\chi(l)B_n(\frac{l+x}{d}). $$
Moreover,
$G_{0,\chi}(x)=0$, and $E_{n, \chi}(x)=\frac{G_{n+1,\chi}(x)}{n+1}.$
\endproclaim
It is not difficult to show that
$$\aligned &\int_{X}\chi(x)e^{(nd+x)t}d\mu(x)-\int_{X}\chi(x)e^{xt}d\mu(x)
=\frac{2\int_{X}\chi(x)e^{xt}d\mu(x)}{\int_{X}e^{ndxt} d\mu(x)}\\
&=\sum_{k=0}^{\infty}\left(2\sum_{l=0}^{dn-1}(-1)^{l-1}\chi(l)l^k\right)\frac{t^k}{k!}.\endaligned\tag17$$
Let $ T_{k, \chi}(n)=\sum_{l=0}^n(-1)^{n-1}\chi(l)l^k.$    From (17), we note that
$$\frac{2\int_{X}\chi(x)e^{xt}d\mu(x)}{\int_{X}e^{ndxt} d\mu(x)}=\sum_{k=0}^{\infty}2T_{k,\chi}(dn-1)\frac{t^k}{k!},
\text{ where $d\in\Bbb N$ with $d\equiv 0$ $(mod  \  2)$, $n\in \Bbb Z_{+}$}, \tag18$$
and
$$\aligned
& \frac{ \int_X \int_X e^{(w_1x_1+w_2x_2)t}\chi(x_1)\chi(x_2) d\mu(x_1)d\mu(x_2)}{\int_X e^{dw_1w_2xt}d\mu(x)}\\
&=\left( \frac{ 2(e^{dw_1w_2t}-1)}{(e^{w_1dt}-1)(e^{w_2dt}-1)}\right) \left(\sum_{a=0}^{d-1}\chi(a)e^{w_1at}(-1)^{a-1}\right)
\left( \sum_{b=0}^{d-1}\chi(b)(-1)^{b-1}e^{w_2bt} \right),\endaligned\tag19$$
where $d\in\Bbb N$ with $d\equiv 0$ $(mod  \  2)$, and $w_1, w_2\in\Bbb N$.

Let $$K(\chi;w_1, w_2|x)=\frac{\int_X \int_X \chi(x_1)\chi(x_2)e^{(w_1x_1+w_2x_2+w_1w_2x)t}d\mu(x_1)d\mu(x_2)}{\int_X e^{d w_1xt}d\mu(x)}.$$
By (18), we see that $K(\chi; w_1, w_2|x)$ is symmetric in $w_1$ and $w_2$.
From(19) and Theorem 4, we obtain the following theorem.
\proclaim{ Theorem 5} For $d\in \Bbb N$ with $d\equiv 0$ $(mod  \  2)$, let $\chi$ be the Dirichlet's character with conductor $d$.
 Then we have
$$\aligned
& \sum_{i=0}^l\binom{l}{i}\frac{d^i}{i+1}\sum_{l=0}^{d-1}(-1)^{l-1}\chi(l)B_{i+1}(\frac{l+w_2 x}{d})T_{l-i, \chi}(dw_1-1) w_1^i w_2^{l-i}\\
&=\sum_{i=0}^{l}\binom{l}{i}\frac{d^i}{i+1}\sum_{l=0}^{d-1}(-1)^{l-1}\chi(l)B_{i+1}(\frac{l+w_1x}{d})T_{l-i,\chi}(dw_2-1)w_2^i w_1^{l-i}.
\endaligned$$
\endproclaim

\enddocument